\documentclass[11pt]{article}
\usepackage{amsfonts,amsmath}
\usepackage{latexsym}
\usepackage{amsthm}
\usepackage{amscd}
\usepackage{commath}
\usepackage{epsfig}
\usepackage{amssymb}
\usepackage{caption}
\usepackage{enumitem}
\usepackage{mathtools}
\usepackage{booktabs}
\usepackage[table]{xcolor}
\usepackage[toc,page]{appendix}
\usepackage[english]{babel}

\usepackage[letterpaper,top=2cm,bottom=2cm,left=2cm,right=2cm,marginparwidth=1.5cm]{geometry}

\usepackage[colorlinks=true, allcolors=blue,bookmarksnumbered]{hyperref}
\hypersetup{
  citecolor = red!80!black
}
\usepackage{color,soul}
\usepackage{float,subcaption}
\usepackage{graphicx}
\usepackage{stackengine}
\usepackage{tikz}
\usetikzlibrary{intersections}
\usepackage[export]{adjustbox}
\usepackage{array}

\usepackage[capitalise]{cleveref}
\usepackage{microtype}

\usepackage{lineno}

\newtheorem{defn}{Definition}[section]
\newtheorem{thm}[defn]{Theorem}

\newtheorem{prop}[defn]{Proposition}
\newtheorem{lem}[defn]{Lemma}

\theoremstyle{remark}

\numberwithin{equation}{section}
\numberwithin{figure}{section}

\newcommand{\bb}{\begin{equation}}
\newcommand{\ee}{\end{equation}}

\newcommand{\elp}{\hspace{-0.7mm}\includegraphics[width=0.35cm]{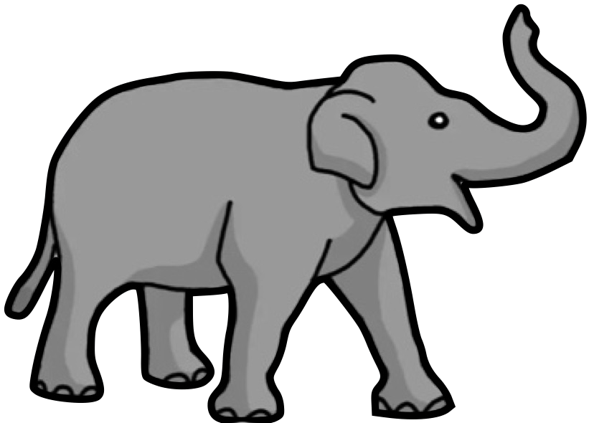}}
\newcommand{\ele}{\hspace{-0.2mm}\includegraphics[width=0.35cm]{elpter}}
\newcommand{\srw}{\includegraphics[width=0.2cm]{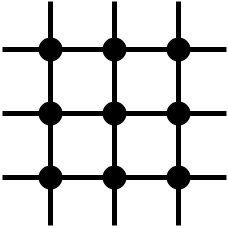}}

\newcolumntype{?}{!{\vrule width 1.5pt }}
\newlength\savedwidth

\newcommand{\tightoverset}[2]{%
  \mathop{#2}\limits^{\vbox to -.5ex{\kern-0.75ex\hbox{$#1$}\vss}}}

\def\rlabel #1 #2{\begin{equation} \label{#1} #2 \end{equation}}

\def\rproof{\begin{proof}}

\def\Qed{\end{proof}}

\makeatletter 
\tikzset{ 
reuse path/.code={\pgfsyssoftpath@setcurrentpath{#1}} 
} 
\tikzset{even odd clip/.code={\pgfseteorule}, 
protect/.code={ 
\clip[overlay,even odd clip,reuse path=#1] 
(current bounding box.south west) rectangle (current bounding box.north east)
; 
}} 
\makeatother 
\usetikzlibrary{3d,perspective,quotes,angles}

\title{\textsc{Recurrence of the plane Elephant random walk}}
\author{
Nicolas \textsc{Curien}\thanks{Universit\'e Paris-Saclay.\hfill  \href{mailto:nicolas.curien@gmail.com}{\texttt{nicolas.curien@gmail.com}}},
Lucile \textsc{Laulin}\thanks{Modal'X -- Université Paris Nanterre \hfill \href{mailto:lucile.laulin@math.cnrs.fr}{\texttt{lucile.laulin@math.cnrs.fr}}}}

\date{}

\begin{document}
\maketitle



\begin{abstract} We give a short proof of the recurrence of the two-dimensional elephant random walk in the diffusive regime. This was recently established by Qin \cite{Qin:2023aa}, but our proof mainly uses very rough comparison with the standard plane random walk. We hope that the method can be useful for other applications. 
\end{abstract}

\section{Introduction} The elephant random walk on $ \mathbb{Z}^d$ has been introduced in dimension 1 by Sch\"utz and Trimper \cite{schutz2004elephants} and is a well-studied discrete process with reinforcement, see \cite{laulin2022elephant} for background and references. Its definition (see \eqref{eq:defelph}) depends on a memory parameter\footnote{The usual definition uses a memory parameter $p \in [0,1]$ which is the probability to reproduce a (uniform) former step of the walk, or to move in one of the $3$ remaining directions with the same probability $(1-p)/3$ so that $\alpha = (4p-1)/3$, see \cite[Eq. (1.4)]{laulin2022elephant}.} $\alpha\in(-\frac{1}{2d-1},1)$ and it exhibits a phase transition going from a diffusive when $\alpha < \alpha_c = \frac{1}{2}$ to a superdiffusive behavior when $\alpha > \alpha_c$. We focus here on the two-dimensional case and establish recurrence of the process in the diffusive regime. 
\begin{thm} 
\label{theo:rec}
In the diffusive regime $ \alpha < \alpha_c = \frac{1}{2}$, the plane elephant random walk is recurrent.
\end{thm}
This has been recently proved by Qin \cite{Qin:2023aa} but our approach is different and much shorter, however it gives less quantitive information and does not directly apply in the critical regime $\alpha = \alpha_c$.
We use a comparaison to the simple random walk which could apply in dimension 1 as well since the simple random walk is also recurrent in that case.

It is worth pointing out that \cite{Qin:2023aa} established that the elephant random walk is always transient when $d \geq 3$, similar to the simple random walk.

\paragraph{Notation.}  We write $  \mathbf{e}_{i}$ the  four directions of $ \mathbb{Z}^2$ for $1 \leq i \leq 4$. We shall write  $ (X_k : k \geq 0)$ for the canonical underlying process starting from $ \mathbf{0}:=(0,0) \in \mathbb{Z}^2$,  we denote its steps by $ \Delta X_k = X_{k+1}-X_{k} \in \{ \mathbf{e}_1, \mathbf{e}_2, \mathbf{e}_3, \mathbf{e}_4\}$ and  we introduce for $1 \leq i \leq 4$  the centered counting direction processes $D^{[X]}_{k}( \mathbf{e}_i)$ defined by
  \begin{eqnarray} D^{[X]}_{k}( \mathbf{e}_i) = \sum_{j=0}^{k-1} \mathbf{1}\{ X_{j+1}-X_{j} = \mathbf{e}_{i}\} - \frac{k}{4}, \qquad \mbox{ in particular notice that }\quad  \sum_{i=1}^4 D^{[X]}_n( \mathbf{e}_i) =0. \label{eq:sum0}  \end{eqnarray} For any stopping time $\theta$, we denote by $X^{(\theta)}$ the shifted process $X^{(\theta)}_k = X_{\theta +k} - X_{\theta}$ for $k \geq 0$. Finally $ \mathcal{F}_n$ is the canonical filtration generated by the first $n$ steps of the walk and we use $X_{[0,n]}$ as a shorthand for $(X_k : 0 \leq k \leq n)$.

\section{Comparison between elephant and simple random walk}
Under the law $ \mathbb{P}_{\srw}$ the underlying process $(X)$ evolves as the standard simple random walk on $ \mathbb{Z}^2$, whereas under $ \mathbb{P}_{ \elp}$, it evolves as the $\alpha$-elephant random walk i.e.\ satisfying for $n \geq 0$
  \begin{eqnarray} \label{eq:defelph} \mathbb{P}_{\elp}(\Delta X_n = \mathbf{e}_i \mid \mathcal{F}_n) = \frac{1}{4} +  \alpha \frac{D^{[X]}_n( \mathbf{e}_i)}{n},  \end{eqnarray}
(where we interpret $0/0=0$ for $n=0$). In particular, under $\mathbb{P}_{\elp}$, the process $( D_k^{[X]}( \mathbf{e}_i) : 1 \leq i \leq 4, k \geq 0)$ is Markov and evolves as an urn process with four colors, which was crucially used in \cite{bertenghi2022functional} to establish the phase transition diffusive/superdiffusive.   The local evolution of the elephant random walk (for large times) ressembles that of the simple random walk and this is quantified in the following propsition:
\begin{prop}[Markov contiguity] \label{lem:contiguity} For any $  \varepsilon>0$ and any $A >0$, there exist $ c_{ \varepsilon, A}>0$  and a sequence of events $ E_n$ satisfying $ \liminf_{n \to \infty} \mathbb{P}_{\srw}(X_{[0,n]} \in E_n) \geq 1- \varepsilon$  such that for any  measurable function $f$, 
$$  \mathbb{E}_{\ele}\left[ f\big( X^{(n)}_{[0,n]}\big)  \mathbf{1}_{ X^{(n)}_{[0,n]} \in E_n} \left| \mathcal{F}_n \mbox{ and } \frac{|D^{[X]}_n( \mathbf{e}_i)|}{ \sqrt{n}}  \leq A \mbox{ for all } 1\leq i\leq4    \right]\right. \geq c_{\varepsilon, A} \cdot  \mathbb{E}_{\srw}\left[ f\big( X_{[0,n]}) \mathbf{1}_{X_{[0,n]} \in E_n}\right].$$
\end{prop}
\proof In the event considered in the conditioning, we have  $ |D_n^{[X]}( \mathbf{e}_i)| \leq A \sqrt{n}$ for all $1 \leq i \leq 4$. By \eqref{eq:defelph}, the Radon--Nikodym derivative of $(X_{n+k}-X_n : 0 \leq k \leq n)$ under $ \mathbb{P}_{\elp}$ with respect to $ \mathbb{P}_{\srw}$ is given by 
 \begin{eqnarray} \label{eq:RND} \mathrm{RND}_n:= \prod_{k=n}^{2n-1} \left( 1 + 4\alpha \frac{D_{k}^{[X]}(\Delta X_k)}{k} \right) = \prod_{k=0}^{n-1} \left( 1 + 4\alpha \frac{D_n^{[X]} ( \Delta X^{(n)}_k) + D_{k}^{[X^{(n)}]}( \Delta X^{(n)}_k)}{n+k} \right).  \end{eqnarray}
By Donsker's invariance principle, we can find a constant $A_ \varepsilon$ such that the event  $$G_n = \left\{ \max_i \sup_{0 \leq k \leq n} |D_k^{[X^{(n)}]}( \mathbf{e}_i)| \leq A_{ \varepsilon} \sqrt{n} \right\}$$ has probability at least $1 - \varepsilon$ under $ \mathbb{P}_{\srw}$. On this event (and conditionally on the event of the statement of the proposition), the counting directions processes $D_{[n,2n]}^{[X]}( \mathbf{e}_i)$ are in absolute value  bounded by $ (A + A_ \varepsilon) \sqrt{n}$. In particular, using $ \log (1 + \alpha x) \geq \alpha x - (\alpha x)^{2}$ for small $|x|$, we deduce that on this event, for $n$ large enough,  the Radon-Nikodym derivative in \eqref{eq:RND} is lower bounded by 
$$ \mathrm{RND}_n \mathbf{1}_{G_n} \geq \exp\left( 4\alpha M_n - (4\alpha)^{2} (A+ A_ \varepsilon)^2\right)\mathbf{1}_{G_n}\quad \mbox{where} \quad M_j = \sum_{k=0}^{j-1} \frac{D_{n+k}^{[X]}(\Delta X_{n+k})}{n+k}.$$ Using \eqref{eq:sum0}, it is trivial to check that $(M_j : 0 \leq j \leq n)$ is a $ (\mathcal{F}_{n+\cdot})$- martingale  with quadratic variation  
$$ \mathbb{E}[M_{j+1}^2-M_j^2\mid \mathcal{F}_{n+j}] = \mathbb{E}[(M_{j+1}-M_j)^2\mid \mathcal{F}_{n+j}] = \frac{1}{4} \frac{\sum_{i=1}^4\left(D_{n+j}^{[X]}( \mathbf{e}_i)\right)^2}{(n+j)^2} \underset{ \mbox{on } G_n}{\leq} \frac{(A+ A_ \varepsilon)^2}{n}.$$
Consider then the $( \mathcal{F}_{n+ \cdot})$-stopping time $$\vartheta = \inf\left\{ j \geq 0: \left\{ \max_i \sup_{0 \leq k \leq j} |D_k^{[X^{(n)}]}( \mathbf{e}_i)| \leq A_{ \varepsilon} \sqrt{n} \right\}\right\}.$$ By the penultimate display, the stopped martingale $(M_{k \wedge \vartheta} : 0 \leq j \leq n)$ has quadratic variation bounded above by $n \times \frac{(A+ A_ \varepsilon)^2}{n}$ and it follows that $  \mathbb{E}[M_n^2 \mathbf{1}_{G_n}]\leq  \mathbb{E}[M_{n\wedge \vartheta}^2] \leq (A+A_{  \varepsilon})^2$. 

In particular, thanks to Markov inequality, for any $ \varepsilon>0$, the event $H_{n} = \{|M_n| \mathbf{1}_{G_n}  < \frac{(A+ A_{ \varepsilon})}{ \sqrt{ \varepsilon}}\}$ is of probability at least $1- \varepsilon$. Gathering up the pieces, on the event $E_n = G_n \cap H_n$ which is of $ \mathbb{P}_{\srw}$ measure at least $1- 2 \varepsilon$, the Radon-Nikodym derivative of the elephant  w.r.t.\ the simple random walk is at least $ \mathrm{e}^{-4\alpha \frac{(A+ A_{ \varepsilon})}{\sqrt{ \varepsilon}}- (4\alpha)^{2} ( A + A_ \varepsilon)^2} =: c_{ \varepsilon,A}$. \endproof

\begin{figure}[!h]
 \begin{center}
 \includegraphics[width=7cm]{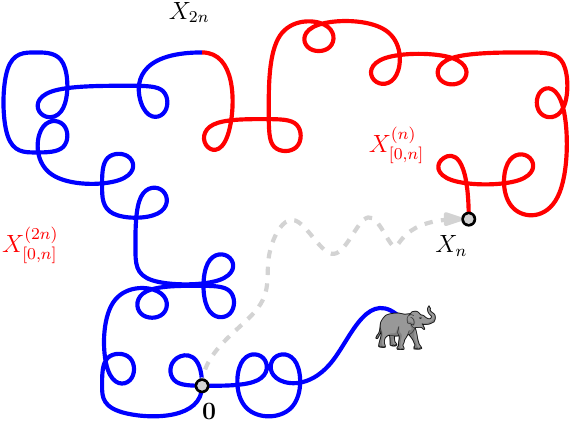}
 \caption{Illustration of the proof of Proposition \ref{prop:retour}. Conditionally on $ \mathcal{F}_n$ and on the fact that the counting directions processes are controlled at time $n$, the blue and red parts are independent on events of large probability. This is sufficient to imply a lower bound on the probability of return to $ \mathbf{0}$. }
 \end{center}
 \end{figure}

 It is classical that in the plane, the simple random walk started from $x \in \mathbb{Z}^2$  with $\|x\| \approx  \sqrt{n}$ has a probability  of order $ \log^{-1} n$ to visit $(0,0)$ within $n$ steps. Our weak bound (Proposition  \ref{lem:contiguity}) is sufficient to imply the same kind of estimate for the elephant random walk:

\begin{prop}  \label{prop:retour}For any $ A>0$ there exists $c_{ A}>0$ such that 
$$\mathbb{E}_{\elp}\left[ \exists \   \frac{5}{2}n\leq k \leq 3n : X_k= 0 \left| \mathcal{F}_n \mbox{ and } \frac{|D^{[X]}_n( \mathbf{e}_i)|}{ \sqrt{n}} \leq A\mbox{ for all } 1\leq i\leq4\right]\right.\geq \frac{c_{A}}{\log n}.$$
\end{prop}

\proof Let us denote $ \mathbf{x}_n = X_n$ which is fixed conditionally on $ \mathcal{F}_n$. Using  Proposition \ref{lem:contiguity} twice, for any positive functions $f$ and $g$ and any $A,A'>0$ and any $ \varepsilon>0$, we can find two sequences of events $E_n$ and $E_n'$ and constants $c_{ \varepsilon,A}$ and $ c_{ \varepsilon,A'}$  such that 
  \begin{eqnarray*}  && \mathbb{E}_{\ele}\left[ f(X^{(2n)}_{[0,n]}) g(X^{(n)}_{[0,n]}) \mid \mathcal{F}_n\right]\\  
  & \geq & \mathbb{E}_{\ele}\left[ f(X^{(2n)}_{[0,n]}) g(X^{(n)}_{[0,n]}) \mathbf{1}_{  \frac{\|D^{[X]}_{2n}( \mathbf{e}_i)\|}{ \sqrt{n}} \leq A',\ \forall 1\leq i\leq4}\Bigg| \mathcal{F}_n\right] \\
  &\geq& c_{ \varepsilon,A'} \cdot \mathbb{E}_{\srw}\left[ f(X_{[0,n]}) \mathbf{1}_{X_{[0,n]} \in E_n'}\right]  \cdot \mathbb{E}_{\ele}\left.\left[ \mathbf{1}_{  \frac{\|D^{[X]}_{2n}( \mathbf{e}_i)\|}{ \sqrt{n}} \leq A',\ \forall 1\leq i\leq4} g(X^{(n)}_{[0,n]}) \right| \mathcal{F}_n\right]\\
  & \geq & c_{ \varepsilon,A'} \cdot \mathbb{E}_{\srw}\left[ f(X_{[0,n]}) \mathbf{1}_{X_{[0,n]} \in E_n'}\right] 
\cdot c_{ \varepsilon,A} \cdot \mathbb{E}_{\srw}\left[ g(X_{[0,n]}) \mathbf{1}_{X_{[0,n]} \in E_n} \mathbf{1}_{\frac{\|D^{[X]}_{n}( \mathbf{e}_i)\|}{ \sqrt{n}} \leq A'-A,\ \forall 1\leq i\leq4}\right] \mathbf{1}_{\frac{\|D^{[X]}_{n}( \mathbf{e}_i)\|}{ \sqrt{n}} \leq A,\ \forall 1\leq i\leq4}.  \end{eqnarray*} 
Up to increasing $A'$ we may suppose that the event $H_n = E_n \cap E'_n \cap \{ \frac{\|D^{[X]}_n( \mathbf{e}_i)\|}{ \sqrt{n}} \leq A'-A,\ \forall 1\leq i\leq4\}$ has probability at least $1- 3 \varepsilon$ and particularizing the inequality above, we deduce that for some constant $ \tilde{c}_{ \varepsilon,A}>0$ the probability in the proposition  is lower bounded by 
$$  \tilde{c}_{ \varepsilon, A} \cdot  \mathbb{P}_{\srw}\left(\exists \   \frac{3}{2}n\leq k \leq 2n : X_k= - \mathbf{x}_n \mbox{ and } \left.\begin{array}{l} X^{(0)}_{[0,n]} \in H_n\\ X^{(n)}_{[0,n]} \in H_n \end{array}\right. \right),$$
so that we can apply the following lemma to conclude.

\begin{lem} For any $A >0$, there exists $ \varepsilon>0$ and $ \delta_{ A}>0$ so that  if $ \mathbf{x}_n \in \mathbb{Z}^2$ is such that $\| \mathbf{x}_n\| \leq A \sqrt{n}$ and if $E_n$ is a sequence of events such that $ \mathbb{P}_{\srw}(X_{[0,n]} \in E_n) \geq 1- \varepsilon$ then we have 
$$ \mathbb{P}_{\srw}\left( \exists  \frac{3}{2}n \leq k \leq 2n : X_{k} = - \mathbf{x}_n \mbox{ and } \left.\begin{array}{l} X^{(0)}_{[0,n]} \in E_n\\ X^{(n)}_{[0,n]} \in E_n \end{array} \right. \right) \geq \frac{\delta_A}{ \log n}.$$
\end{lem}
\proof We use a second-moment method on the random variable  
$$ \mathcal{N}_{ \mathbf{x}_n}^{E_n} := \# \left\{ \frac{3}{2} n \leq k \leq 2n : X_k = - \mathbf{x}_n\right\} \mathbf{1}_{X^{(0)}_{[0,n]} \in E_n} \mathbf{1}_{X^{(n)}_{[0,n]} \in E_n}.$$ We denote by $p_k^{E_n}(y) = \mathbb{E}_{\srw}[ \mathbf{1}_{X_k = y} \mathbf{1}_{X_{[0,n]} \in E_n}]$ and $p_k(y) = \mathbb{P}_{\srw}(X_k = y)$ for the heat kernels. By the standard local limit theorem (or just Stirling approximation on the binomial coefficients) there exists $C>0$ such that $p_k(y) \leq \frac{C}{k}$ for all $k \geq 1$ and $y \in \mathbb{Z}^2$. First, by lifting the restrictions on $E_n$ we have 
 \begin{eqnarray*} \mathbb{E}_{\srw}\left[ \left(\mathcal{N}_{ \mathbf{x}_n}^{E_n}\right)^2\right]  \leq  \mathbb{E}_{\srw}\left[ \left(\sum_{k=3/2n}^{2n} \mathbf{1}_{X_k =  -\mathbf{x}_n}\right)^2\right] &\leq& 2\sum_{ \frac{3}{2}n \leq k \leq k' \leq 2n} p_k( -\mathbf{x}_n) p_{k'-k}( \mathbf{0}) \\ &\leq&  2\sum_{ \frac{3}{2}n \leq k \leq k' \leq 2n}  \frac{C}{n} \frac{C}{k'-k} \leq {C'} \log (n),
 \end{eqnarray*}
for some $C' >0$ (independent of $n$). To evaluate the first moment,  introduce the (truncated) Green functions $g^{E_n}(y) =  \sum_{k=n/2}^{n} p_k^{E_n}(y)$ and similarly $g(y)= \sum_{k=n/2}^{n} p_k(y)$. In particular, since $ \mathbb{P}_{\srw}(E_n) \geq 1- \varepsilon$ we have $\|p-p^{E_n}\|_1 := \sum_y p(y)-p^{E_n}(y) \leq \varepsilon$ and similarly and $\|g-g^{E_n}\|_1=\sum_y g(y)-g^{E_n}(y) \leq \varepsilon n$. Recalling that $ \frac{C}{n} \geq p^{E_n}_n(y) \geq p_n(y) $ and $2C \geq g^{E_n}(y)\geq g(y)$, we have 
 \begin{eqnarray*} \mathbb{E}_{\srw}[ \mathcal{N}_{ \mathbf{x}_n}^{E_n}] = \sum_{y \in \mathbb{Z}^2}p_n^{E_n}(y) g^{E_n}(-y- \mathbf{x}_n) & =& 
 \sum_{y \in \mathbb{Z}^2}\left( \begin{array}{l}  p_n^{E_n}(y) g^{E_n}(-y- \mathbf{x}_n) - p^{E_n}_n(y)g(-y- \mathbf{x}_n)\\
- p_n(y)g(-y- \mathbf{x}_n) + p_n^{E_n}(y)g(-y- \mathbf{x}_n) \\
  +p_n(y)g(-y- \mathbf{x}_n)
  \end{array}\right)\\
& \geq & \sum_{y}p_n(y)g(-y- \mathbf{x}_n) - \|p_n^{E_n}\|_{\infty}\|g-g^{E_n}\|_1 - \|g\|_\infty \|p_n-p_n^{E_n}\|_1\\  &\geq& \sum_{y}p_n(y)g(-y- \mathbf{x}_n) - 3C^2 \varepsilon.  \end{eqnarray*}
However, since $ \|\mathbf{x}_n\| \leq A \sqrt{n}$, the local limit theorem implies that $\sum_{y}p_n(y)g(-y- \mathbf{x}_n)> c_A$ for some $ c_A>0$ independently of $n$ and so one can choose $ \varepsilon>0$ small enough so that if $ \mathbb{P}_{\srw}(E_n) \geq 1- \varepsilon$ then we have  $\mathbb{E}_{\srw}\left[ \mathcal{N}_{ \mathbf{x}_n}^{E_n}\right]>c_A/2$. We conclude by the second moment method that
$$ \mathbb{P}_{ \srw}( \mathcal{N}_{ \mathbf{x}_n}^{E_n}>0) \geq \mathbb{E}_{\srw}\left[ \mathcal{N}_{ \mathbf{x}_n}^{E_n}\right]/\mathbb{E}_{\srw}\left[ \left(\mathcal{N}_{ \mathbf{x}_n}^{E_n}\right)^2\right] \geq \frac{c_A}{2 C'\  \log n}.$$ \endproof

\section{Proof of Theorem \ref{theo:rec}}
\proof[Proof of Theorem \ref{theo:rec}]
 Let us denote  $ \mathcal{P}_{3^j} =  \mathbb{P}_{\elp}( \exists \  3^{j}\leq k \leq 3^{j+1}, X_k= \mathbf{0} \mid \mathcal{F}_{3^j})$. When $\alpha < \alpha_{c}$, i.e. the diffusive regime,  Bertenghi  \cite[Theorem 4.2]{bertenghi2022functional} showed  that under $ \mathbb{P}_{\ele}$ we have 
$$ \left( \frac{D_{n}^{X}( \mathbf{e}_{i})}{ \sqrt{n}}\right)_{1 \leq i \leq 4} \xrightarrow[n\to\infty]{(d)} ( \mathcal{X}_{i})_{1 \leq i \leq 4},$$ for some random variable $ \mathcal{X}$ (whose distribution is irrelevant for our purposes). Together with our Proposition \ref{prop:retour}, this shows that in the diffusive regime, for any $ \varepsilon>0$ there exists $ \delta>0$ such that for large $j$'s we have
  \begin{eqnarray} \label{eq:atom0} \mathbb{P}_{\ele}( j \cdot \mathcal{P}_{3^j} > \delta) \geq 1-\varepsilon.  \end{eqnarray} Notice that the variables $ \mathcal{P}_{3^j}$ are not independent, but Jeulin's lemma \cite[Proposition 3.2]{matsumoto2011zero} gives 
  \begin{eqnarray} \label{eq:jeulin} \sum_{k=1}^{\infty}  \mathcal{P}_{3^{k}} = \infty, \qquad  \mathbb{P}_{\ele}-a.s.  \end{eqnarray}
  To be honest we rather use the proof that the lemma itself, and since the argument is short let us reproduce it here: Suppose by contradiction that there exists $ \varepsilon,M>0$ so that the event $A= \{\sum_{k=1}^{\infty}  \mathcal{P}_{3^{k}} < M\}$ has probability at least $ \varepsilon>0$. Using \eqref{eq:atom0} we take $\delta >0$  so that  $\mathbb{P}_{\ele}( j \cdot \mathcal{P}_{3^j} > \delta) \geq 1- \frac{\varepsilon}{2}$ and write
  $$ M \geq \mathbb{E}\left[   \mathbf{1}_{A} \sum_{j \geq 1} \mathcal{P}_{3^{k} }\right] \geq \sum_{j \geq 1} \frac{\delta}{j} \cdot \underbrace{\mathbb{P}_{\elp}\left(A \cap \left\{\mathcal{P}_{3^j} > \frac{\delta}{j} \right\}\right)}_{ \geq \varepsilon - (1 - (1- \frac{ \varepsilon}{2})) = \varepsilon/2} = \infty,$$  which is a contradiction. Given \eqref{eq:jeulin}, the conditional Borel-Cantelli lemma  (\cite[Theorem 4.3.4]{durrett2019probability}) then implies that the events $\{\exists \  3^{j}\leq k \leq 3^{j+1} : X_k= \mathbf{0}\}$ happen for infinitely many $j$'s with probability one, implying recurrence of the process.

\bigskip 
\noindent \textbf{Acknowledgments.} We thank Jean Bertoin and Thomas Budzinski for  discussions about the plane elephant random walk. We are also indebted to Alice Contat for providing us with the pictures $\srw$ and $\elp$. We are grateful to Zheng Fang for a very careful reading of the paper.

\bibliographystyle{siam}
\bibliography{bib-rec}

\end{document}